\documentclass[10pt]{article}

\usepackage[utf8]{inputenc}
\usepackage[T1]{fontenc}
\usepackage[english]{babel}
\usepackage{lmodern}

\usepackage{amsmath,amssymb,bm}

\usepackage{graphicx}
\usepackage{booktabs}

\usepackage[margin=1.5in]{geometry}
\makeatletter
\renewcommand\footnoterule{%
  \vspace{0.7em}%
  \hrule width 0.4\columnwidth height 0.4pt
  \kern 0.3em%
}
\makeatother

\usepackage[hidelinks]{hyperref}

\begin{document}

\title{\textbf{Monetisation and Cost Function Scalarisation Cannot Guarantee Valid Multi-Criteria Decision-Making}}

\author{
A.R.M. Wolfert\thanks{Corresponding author: A.R.M. (Rogier) Wolfert: \texttt{a.r.m.wolfert@tudelft.nl}}\\[0.4em]
\small Department of Algorithmics, Faculty of Mathematics and Computer Science\\[-0.3em]
\small Delft University of Technology, The Netherlands
}

\date{}
\maketitle
\vspace{-2.0em}

\begin{abstract}

\noindent
Monetary quantities are not preferences. Financial cost is an object-related system performance, whereas value-for-money expresses an actor's subjective evaluation of that cost. Monetisation embeds subjective trade-off information: its conversion coefficients express these valuations and determine how heterogeneous system performances compensate for one another. Cost-benefit analysis converts performances into monetary equivalents, but the results remain quantities of money. Even net present value is a discounted summation of monetary quantities, not a measurement of the preference value that actors attribute to them.

\noindent
This creates a monetisation paradox. Cost initially enters a multi-criteria problem as one performance among several, but money is subsequently elevated into the universal representation of all other performances, thereby creating representational dominance. The original multi-criteria problem is consequently reformulated as a single-criterion problem in which all performances influence the outcome only through their imposed monetary equivalents. Cost-function scalarisation does not resolve this problem by replacing money with points, normalised scores, losses, or penalties. Such transformations and their weights determine preference directions, intensities, and compensatory trade-offs and therefore constitute implicit preference models, although they are commonly treated as technical scalar constructions rather than explicitly recognised, measured, and validated as preference models.

\noindent
A mathematically well-defined scalar optimum is consequently not necessarily a decision-valid multi-criteria outcome. Numerical commensurability does not establish preference commensurability or a unique preference-consistent ranking. Decision validity requires heterogeneous system performances to be transformed into actor preferences represented in a unified preference domain. Preference is relational and contextual rather than an absolute scalar quantity, but it can be represented numerically on an affine preference scale. Preference Function Modelling and the Linear Preference Space provide the common reference and commensurability required for valid preference aggregation. Cost-function scalarisation forms a weighted sum of transformed performance measures, whereas decision-valid MCDM aggregates validly modelled and commensurate actor preferences concerning those performances. The final decision criterion should therefore maximise aggregated actor preference over the feasible and acceptable design--decision space rather than monetised quantities or scalarised performance scores.

\end{abstract}

\vspace*{1.0em}


\vspace{0.75em}

\noindent\textbf{Article history:}\\
Compiled November~03,~2026

\newpage

\section*{1. Introduction}

Multi-criteria decision-making (MCDM), multi-criteria decision analysis (MCDA), and multi-objective optimisation (MOO) address problems in which alternatives or system designs are characterised by several, often conflicting, performances. Established MCDA provides structured approaches for comparing alternatives across criteria, while value-focused thinking emphasises the values and objectives that should guide the decision \cite{belton2002,keeney1992}.

A design \(x\) generates a system-performance vector

\begin{equation}
\label{eq:performance-vector}
\mathbf{F}(x)
=
\left[
f_1(x),\ldots,f_I(x)
\right].
\end{equation}

For example, \(\mathbf{F}(x)\) may comprise
\(\left[C(x),T(x),E(x),S(x),Q(x),\ldots\right]\), where \(C(x)\) denotes cost, \(T(x)\) duration, \(E(x)\) environmental impact, \(S(x)\) safety, and \(Q(x)\) quality or reliability. These functions describe what the system performs, but they do not by themselves establish how desirable those performances are to the relevant actors. System performances are object-related consequences, whereas preferences are subject-related evaluations of those consequences.

Several established disciplines provide methods for modelling, exploring, and comparing such alternatives. Classical MOO studies conflicting and potentially incommensurable objectives, Pareto optimality, trade-offs, and the role of the decision-maker \cite{miettinen1999}. Evolutionary MOO uses population-based search to approximate sets of Pareto-optimal alternatives, particularly for complex and nonlinear design spaces \cite{deb2001}. These approaches are powerful for exploring feasible performance spaces, but the generation of a Pareto set does not itself identify which alternative provides the most preferred or decision-valid outcome \cite{ozerwolfert2026}.

Systems engineering adopts a total life-cycle perspective extending from needs and requirements through design, operation, support, phase-out, and disposal \cite{blanchardfabrycky2011}. Life-cycle cost analysis and engineering economy compare alternatives using monetary quantities, economic equivalence, present worth, rate of return, and benefit-cost analysis \cite{fabryckyblanchard1991,thuesenfabrycky2001}. Value engineering similarly evaluates alternatives through their functions, performance, quality, risks, and costs \cite{dellisola1997,park1999}. These methods can provide valid economic and object-level comparisons, but their numerical results do not automatically measure the preference value that actors attribute to the systems or consequences being compared.

Cost-benefit analysis represents the established monetisation route. It converts heterogeneous consequences into monetary equivalents, discounts future costs and benefits, and compares alternatives through a scalar monetary criterion such as net present value \cite{boardman2018,aldy2021}. Monetisation is applied not only to direct expenditures and revenues, but also to environmental, health, safety, and other consequences. Multi-criteria analysis has therefore also been positioned as an alternative for evaluating heterogeneous consequences without requiring their compulsory translation into money \cite{gamper2015}.

Cost-function scalarisation provides a broader route in which performances are transformed into monetary quantities, normalised scores, points, losses, penalties, utilities, or other numerical representations and subsequently combined into one scalar objective \cite{marler2004,ehrgott2005,martins2022}. Monetisation and cost-function scalarisation are not equivalent. Monetisation expresses heterogeneous performances in monetary units, whereas non-monetary scalarisation uses composite numerical scores. Nevertheless, both embed implicit preference and trade-off information. Monetary conversion coefficients determine how heterogeneous performances compensate for one another, while score transformations, penalties, aggregation rules, and weights determine analogous trade-offs numerically.

Subjective information may therefore be embedded in an apparently objective scalar criterion without being explicitly recognised, measured, or validated as actor preference. This creates two related validity problems. First, monetisation begins with financial cost as one system performance among several but elevates money into the universal representation of all performances. Second, cost-function scalarisation replaces heterogeneous performances with a single constructed numerical objective whose effective trade-offs depend jointly on score construction and weighting. Neither monetary commensurability nor a common numerical score range establishes preference commensurability. These arguments are developed for monetisation in Section~2 and for cost-function scalarisation in Section~3.

MCDM should therefore not be reduced to the mechanical application of a numerical method to achieve computational optimality. French \cite{french2023} places formal multi-criteria modelling within a broader and iterative decision process involving sense-making and modelling, analysis and exploration, and interpretation and implementation. This perspective is directly relevant to the distinction between computational optimality and decision validity, particularly when monetary conversion coefficients, scores, penalties, transformations, and weights are introduced as technical inputs although they contain substantive and actor-dependent trade-off information.

A physical system has performances but no preferences or decision purpose of its own. Goal orientation arises when actors evaluate system performances:

\begin{equation}
\label{eq:actor-preference}
P_{k,i}\bigl(f_i(x)\bigr),
\end{equation}

where \(P_{k,i}\) is the performance-to-preference function of actor \(k\) concerning performance \(i\). This distinction is central to preference measurement and preference-based decision theory \cite{barzilai2005,barzilai2010,barzilai2022}. Cost is therefore one system performance that must itself be evaluated through actor preference, rather than automatically elevated into the universal measure of all performances.

The problem is not the use of a numerical decision criterion as such. The decisive question is what that criterion validly represents. A numerical criterion may represent a physical performance, a monetary consequence, a weighted point score, or an actor preference. These quantities are not interchangeable. Numerical precision and computational optimality do not by themselves establish preference validity or decision validity.

In this article, \textbf{decision validity} denotes the ability of an MCDM method to produce a preference-valid, consistent, and unique decision outcome. Decision validity requires actor preferences to be validly modelled, represented in a common preference domain with a stable affine reference origin and a commensurate preference unit, and subsequently aggregated consistently. Equivalent preference information must consequently produce the same ranking and selected alternative irrespective of arbitrary numerical representations.

A decision-valid outcome is therefore a unique best-fit-for-common-purpose solution determined by validly modelled and consistently aggregated actor preferences, rather than by chosen physical units, monetary conversion coefficients, score ranges, normalisation rules, or arbitrary scale-weight combinations \cite{barzilai2005,barzilai2010,wolfert2026unique}. This requires an explicit transformation from system performance to actor preference, a distinction between technical feasibility and actor acceptability, and valid preference aggregation over the feasible and acceptable solution space \cite{wolfert2023,wolfert2026unique}.

A mathematically well-defined scalar optimum is therefore not necessarily a decision-valid multi-criteria outcome. Attempts to validate an unvalidated scalarisation retrospectively through Pareto efficiency or sensitivity analysis create a two-stage decision-validity failure. First, monetisation or cost-function scalarisation imposes trade-offs without establishing that they validly represent actor preferences. Second, performance-space ordering or parameter variation is used to support the resulting optimum, although neither can supply or validate the missing preference foundation. The detailed two-stage failure is examined in Section~3.2.

This article consequently examines whether these two scalar routes can guarantee decision validity. Monetisation produces monetary quantities whose conversion coefficients contain implicit actor-dependent trade-offs, while cost-function scalarisation constructs numerical criteria whose transformations and weights determine analogous trade-offs. Neither route automatically establishes valid preference measurement or aggregation, and its computational optimum cannot therefore be assumed to constitute a decision-valid multi-criteria outcome. The central question is thus not whether a scalar-valued criterion may be used, but whether that criterion represents validly measured and consistently aggregated actor preference. Performance scalarisation constructs a scalar objective directly from system performances. Preference aggregation constructs a numerical decision criterion only after actor preferences have been explicitly modelled, transformed into a common preference domain, made commensurate, and validly aggregated.

Section~2 examines why monetisation cannot guarantee decision validity. Section~3 develops the corresponding argument for cost-function scalarisation. Section~4 presents the performance-to-preference transformation and the valid preference-aggregation alternative.

\section*{2. Why Monetisation Is Decision-Invalid}

Monetisation creates a fundamental paradox in multi-objective optimisation. In the original multi-criteria formulation, financial cost has the same logical status as duration, emissions, safety, quality, reliability, and other system performances:

\begin{equation}
\label{eq:performance-vector-monetisation}
\mathbf{F}(x)
=
\left[
C(x),T(x),E(x),S(x),Q(x),\ldots
\right].
\end{equation}

Cost is therefore one component of the system-performance vector, not the natural measure of its other components. It describes one object-related consequence of a design, just as duration, emissions, safety, and quality describe other consequences. A monetised formulation changes this relationship:

\begin{equation}
\label{eq:monetised-function}
M(x)
=
C(x)
+\lambda_T T(x)
+\lambda_E E(x)
+\lambda_S S(x)
+\lambda_Q Q(x).
\end{equation}

The coefficients \(\lambda_i\) convert non-monetary performances into monetary equivalents. Duration becomes the cost of time, emissions the cost of emissions, safety the cost of safety consequences, and quality the cost of insufficient quality. A performance can therefore influence the scalar result only after being translated into monetary terms.

Monetisation consequently begins with cost as one performance among several but ends with money as the representation of every performance. Cost itself does not become every performance; rather, its monetary dimension is elevated into the universal representation through which all performances must be expressed. This creates \textbf{representational dominance}. It does not imply that the direct financial-cost term \(C(x)\) is numerically dominant, since a monetised safety, environmental, or time-related term may be larger. It means that every performance must be represented in one privileged dimension, namely money.

This hierarchy is representational rather than lexicographic. The performances are not optimised successively according to an explicit priority order. Instead, money is selected as the governing representation, the coefficients \(\lambda_i\) determine the compensatory trade-offs, and the resulting monetary quantities are summed and optimised as a single objective:

\begin{equation}
\label{eq:monetised-optimum}
x_M^*
=
\arg\min_{x\in X} M(x).
\end{equation}

The result is an optimum of the constructed monetary model, not necessarily an optimum of validly measured and aggregated actor preferences. Although the conversion coefficients may contain implicit preference and trade-off information, their outputs remain quantities of money.If \(m\) denotes a monetary quantity, the actor-specific and context-dependent preference attributed to it must be represented separately as \(P_k(m)\neq m\): \(P_k(m)\) represents actor \(k\)'s preference evaluation, whereas \(m\) represents an amount of currency. Monetary arithmetic can therefore aggregate quantities of money, but it does not thereby aggregate actor preferences.

This distinction also applies to cost-benefit analysis. A generic net present value function is

\begin{equation}
\label{eq:npv}
NPV(x)
=
\sum_{t=0}^{T}
\frac{B_t(x)-C_t(x)}
{(1+r)^t},
\end{equation}

where \(B_t(x)\) and \(C_t(x)\) are monetised benefits and costs, \(r\) is the discount rate, and \(T\) is the assessment horizon. Despite the word \emph{value} in its name, NPV is a discounted summation of monetary quantities, not a measurement of the preference value that actors attribute to them. NPV therefore establishes monetary summation, not preference aggregation \cite{boardman2018,barzilai2005,barzilai2010,barzilai2022}.

Monetisation establishes \textbf{accounting commensurability} because all terms are expressed in the same currency unit. It does not establish \textbf{preference commensurability}. Equal monetary quantities need not represent equal preference differences across consequences or actors:

\begin{equation}
\begin{aligned}
P_k(\text{EUR }1\text{ environmental consequence})
&\neq
P_k(\text{EUR }1\text{ safety consequence}),
\\
P_A(\text{EUR }1)
&\neq
P_B(\text{EUR }1).
\end{aligned}
\end{equation}

Likewise, the preference value of a monetary sum is not generally equal to the sum of the preference values of its components:

\begin{equation}
\label{eq:nonadditive-value}
P_k(a+b)
\neq
P_k(a)+P_k(b).
\end{equation}

Equality would require additional assumptions, including linear preference for money over the complete decision-relevant range. \textsc{Note:} Although both \(\text{EUR }11\text{m}-\text{EUR }10\text{m}\) and \(\text{EUR }101\text{m}-\text{EUR }100\text{m}\) equal \(\text{EUR }1\text{m}\), the corresponding preference differences need not be equal. Identical monetary differences can therefore have different preference meanings depending on the actor, context, and location on the monetary scale.

Cost-benefit analysis can validly sum monetary amounts as monetary quantities; whether those amounts constitute meaningful and commensurate preference differences is a separate measurement question.

Monetisation is nevertheless not preference-free. A coefficient such as

\begin{equation}
\lambda_E
=
\text{EUR }100/\text{tonne CO}_2
\end{equation}

does not arise solely from the physical properties of emissions. It may depend on the valuation method, affected population, willingness and ability to pay, time horizon, discount rate, risk assumptions, institutional context, and distribution of consequences. Cost-benefit analysis may use market prices, shadow prices, revealed preferences, stated preferences, contingent valuation, or values of avoided consequences \cite{boardman2018,aldy2021}. The transformation is therefore more accurately represented as \(\text{system performance}\longrightarrow\text{human or institutional valuation}\longrightarrow\text{monetary equivalent}\). The conversion coefficients encode assumptions about relative importance and exchangeability. The optimiser does not discover these trade-offs; it executes those encoded in the coefficients \(\lambda_i\).

\begin{quote}
\textbf{A monetised performance is a monetary quantity containing implicit trade-off assumptions, but it is not thereby a validated preference measurement.}
\end{quote}

Unit and currency invariance are necessary for the internal consistency of a monetary model but remain insufficient for preference validity. If \(1\text{ EUR}=r\text{ USD}\), a consistently converted monetary function satisfies

\begin{equation}
M_{\$}(x)
=
rM_{\text{EUR}}(x),
\qquad
r>0,
\end{equation}

and therefore

\begin{equation}
\arg\min_x M_{\text{EUR}}(x)
=
\arg\min_x M_{\$}(x).
\end{equation}

This establishes that different numerical amounts can represent the same monetary consequence. It does not establish that the consequence itself measures preference.

Actual financial cost remains an important special case. If every term represents an actual expenditure or revenue, the monetary quantities may legitimately be added:

\begin{equation}
\label{eq:total-cost}
C_{\mathrm{total}}(x)
=
\sum_i p_i q_i(x).
\end{equation}

Minimising total financial cost is then a valid object-level optimisation problem:

\begin{equation}
x_C^*
=
\arg\min_{x\in X}
C_{\mathrm{total}}(x).
\end{equation}

This identifies the feasible system with the lowest total financial cost. It does not necessarily identify the feasible and acceptable system most preferred by the relevant actors:

\begin{equation}
C_{\mathrm{total}}(x)
\neq
P_{k,C}\bigl(C_{\mathrm{total}}(x)\bigr).
\end{equation}

A minimum-cost solution constitutes the decision optimum only when minimum financial cost has explicitly been established as the complete decision purpose. Otherwise, it remains an object-performance optimum.

The distinction becomes particularly important in group decision-making because actors may evaluate the same performance differently:

\begin{equation}
P_{A,i}(f_i)
\neq
P_{B,i}(f_i).
\end{equation}

Actors may differ in preference direction, intensity, thresholds, target values, acceptable intervals, sensitivity to gains and losses, and relative importance. An operations actor may prefer a longer project duration to obtain execution margin, while a commercial actor may prefer a shorter duration to obtain earlier revenue. A single monetary cost per day replaces both preference functions with one conversion rule \(\lambda_TT\).

Preferences may also exhibit saturation, tipping points, asymmetric penalties, acceptable ranges, and interior optima. For duration, for example,

\begin{equation}
P_T(T)
=
\begin{cases}
\text{increasing}, & T<T^*,\\
\text{maximum},    & T=T^*,\\
\text{decreasing}, & T>T^*.
\end{cases}
\end{equation}

A constant monetary coefficient imposes a fixed marginal conversion and generally a fixed direction of improvement. It cannot preserve such preference structures unless the relevant preference is genuinely linear over the complete domain.

The suppression of actor-specific preference structures also explains why monetary compensation does not establish acceptability. A monetised sum may imply

\begin{equation}
+\text{EUR }1\text{m benefit to actor A}
-
\text{EUR }1\text{m impact on actor B}
=
0.
\end{equation}

This arithmetic cancellation does not establish that actor B accepts the outcome, that both actors experience equal preference differences, that the consequences are legitimately exchangeable, or that the result serves the common purpose. Technical feasibility and actor acceptability are distinct conditions \cite{Wolfert2026a}. Feasibility may be written as

\begin{equation}
g_f\bigl(\mathbf{F}(x)\bigr)
\leq
0,
\end{equation}

whereas actor-specific acceptability requires

\begin{equation}
P_{k,i}\bigl(f_i(x)\bigr)
\geq
\underline{P}_{k,i}.
\end{equation}

The relevant design--decision space is therefore 

\begin{equation}
\label{eq:feasible-acceptable-space}
\mathcal{S}_{f,a}
=
\left\{
x\in X:
g_f\bigl(\mathbf{F}(x)\bigr)\leq0
\ \land\
P_{k,i}\bigl(f_i(x)\bigr)\geq\underline{P}_{k,i}
\quad
\forall\,k,i
\right\}.
\end{equation}

Monetary compensation may conceal unacceptable outcomes by allowing favourable monetary terms to offset an actor's unacceptable performance. A decision-valid multi-criteria method must instead preserve the distinctions between monetary consequence and actor preference, accounting and preference commensurability, and technical feasibility and actor acceptability \cite{Wolfert2026a}. Monetisation can support economic analysis and object-level optimisation, but it cannot by itself guarantee a valid multi-actor decision outcome.

\section*{3. Why Cost-Function Scalarisation Is Decision-Invalid}

Cost-function scalarisation transforms multiple heterogeneous system performances, whether monetised, scored, penalised, normalised, or otherwise transformed, into a single optimisation criterion. It embeds implicit preference assumptions and compensatory trade-offs, yet does not by itself satisfy the decision-validity requirements defined in Section~1. The central question is therefore not whether this scalar criterion can be optimised, but whether it validly represents and aggregates actor preferences.

\subsection*{3.1 Scalar Cost Functions Are Implicit but Unvalidated Preference Models}

In optimisation, the term \emph{cost function} is often used generically for a scalar objective to be minimised:

\begin{equation}
\label{eq:generic-cost-function}
x^*
=
\arg\min_{x\in X}J(x).
\end{equation}

The function \(J(x)\) may represent financial cost, physical mass, energy consumption, prediction error, environmental impact, distance, penalties, utility scores, or a weighted combination of several performances. Every monetised function can therefore be represented as a scalarised cost function, but not every scalarised cost function is monetised. A cost function reduces multiple performances to one optimisation criterion. Its computational optimisability does not establish that it validly represents actor preference.

Replacing euros with points does not create a preference-neutral model. Suppose each performance is transformed into a score \(s_i\bigl(f_i(x)\bigr)\in[0,100]\) and subsequently combined in a weighted cost function:

\begin{equation}
\label{eq:scalarised-cost-function}
J(x)
=
\sum_{i=1}^{I}
w_i s_i\bigl(f_i(x)\bigr).
\end{equation}

Although \(s_i\) may be described as a technical score, normalisation, penalty, loss, utility, fitness, or evaluation function, the function contains substantive preference and trade-off information. It determines the preferred direction of change, the numerical intensity of improvement or deterioration, marginal effects, targets or thresholds, and the representation of gains and losses. The weights subsequently determine how the resulting score differences compensate across performances.

A weighted cost function is therefore an implicit preference model. Its scores are numerical expressions of desirability, and their introduction moves the model from the object-performance domain towards the subject-preference domain. If this transition is not recognised, the preference validity of the resulting scalar objective remains unexamined.

Assigning numerical scores does not, however, establish preference measurement. If one outcome receives \(80\) points and another \(60\), their numerical difference is \(20\). This arithmetic difference does not explain what \(20\) preference points mean empirically or establish that equal score differences represent equal preference differences. The mere numerical representation of preference does not justify addition, multiplication, averaging, distance calculation, or aggregation without an appropriate measurement structure \cite{barzilai2005,barzilai2010}.

Preference representations are admissible under positive affine transformations:

\begin{equation}
\label{eq:affine-transform}
P'(x)
=
aP(x)+b,
\qquad
a>0.
\end{equation}

The numerical zero and unit are therefore reference choices rather than absolute properties of preference. Preference measurement must precede aggregation. A scalar objective is not automatically a preference measure merely because it is expressed numerically.

Normalisation does not resolve this problem. Suppose several performances are normalised such that \(s_i(f_i)\in[0,1]\). The resulting scores share a numerical interval, but this does not imply

\begin{equation}
0.1\text{ cost points}
=
0.1\text{ safety points}
=
0.1\text{ environmental points}
\end{equation}

as preference differences. Normalisation aligns numerical endpoints and ranges; it does not establish a common preference unit. The same numerical range does not establish the same preference scale, and numerical comparability is not preference commensurability.

Moreover, the effective trade-offs in Equation~\eqref{eq:scalarised-cost-function} are determined jointly by score construction and weighting. If a score is transformed as

\begin{equation}
\label{eq:score-rescaling}
s_i'(f_i)
=
a_i s_i(f_i)+b_i,
\qquad
a_i>0,
\end{equation}

then its weighted contribution becomes

\begin{equation}
w_i s_i'(f_i)
=
w_i a_i s_i(f_i)+w_i b_i.
\end{equation}

Preserving the original effective scaling requires the associated weight to be adjusted according to

\begin{equation}
w_i'
=
\frac{w_i}{a_i}.
\end{equation}

A trade-off is therefore not determined by the weight alone. It is jointly determined by the score function and its weight. The same numerical weight may encode different trade-offs under different normalisations, while different combinations of score scales and weights may encode the same intended trade-off. Consequently, different defensible scalarisations can produce different decisions without any underlying actor preference having changed.

For example,

\begin{equation}
J_1(x)
=
0.5s_1\bigl(f_1(x)\bigr)
+
0.5s_2\bigl(f_2(x)\bigr)
\end{equation}

may appear to assign equal importance to both performances. Yet if the scores use different ranges, reference values, or normalisation rules, one term may dominate numerically. Altering these choices can change the selected optimum while retaining the same stated weights.

The optimiser may nevertheless return a mathematically unique result, but this proves only optimality relative to the selected numerical construction. It does not establish that the result uniquely maximises actor preference. A unique numerical optimum does not imply a unique preference optimum, and computational uniqueness does not establish decision validity.

The limitation becomes more fundamental in multi-actor decisions. A single scalar cost function commonly represents the problem as though one preference structure applies to the entire group, whereas actor-specific preferences may differ:

\begin{equation}
P_{A,i}(f_i)
\neq
P_{B,i}(f_i).
\end{equation}

Actors may differ in preferred directions, preference intensities, acceptable ranges, targets, trade-offs, and priorities. Actors may even hold opposing preferences concerning the same system performance. Combining transformed performances into one scalar score before explicitly representing these differences suppresses the actor structure.

A decision-valid group model must instead preserve \(P_{k,i}\bigl(f_i(x)\bigr)\), so that actor \(k\) and performance \(i\) remain explicit until valid preference aggregation takes place.

Cost-function scalarisation therefore creates numerical commensurability, but it does not necessarily establish preference measurement, preference commensurability, or valid multi-actor aggregation. Decision-valid MCDM instead requires explicitly modelled actor preferences to be represented in a common preference domain, made commensurate, and validly aggregated.

\begin{quote}
\textbf{A weighted scalar cost function embeds subjective preference information but does not constitute valid preference aggregation.}
\end{quote}.

\subsection*{3.2 Pareto Analysis Cannot Repair an Invalid Scalarisation}

Once a monetised or weighted cost function has been specified, its conversion coefficients, score transformations, penalties, and weights already define compensatory trade-offs and thereby determine the scalar optimum \(x_J^*=\arg\min_{x\in X}J(x)\).

The original multi-performance problem has then effectively been transformed into a single-objective optimisation problem. Pareto analysis is not required to compute this scalar optimum. It may still support performance-space exploration, generate non-dominated alternatives, or examine sensitivity to changes in coefficients and weights. Classical and evolutionary MOO commonly use Pareto dominance for these purposes \cite{ehrgott2005,HillierLieberman2021,deb2001,martins2022}.

Demonstrating that a scalar optimum is Pareto-efficient does not validate the scalarisation through which the optimum was obtained. Pareto dominance orders performances under assumed directions of improvement; it neither establishes nor validates the preference information embedded in monetary conversion coefficients, scoring functions, penalties, transformations, or weights. When these assumed improvement directions do not preserve the ordering induced by aggregated actor preference, Pareto filtering may even exclude the best-fit-for-common-purpose solution and thereby produce a decision-invalid outcome \cite{ozerwolfert2026}.

When an unvalidated scalarisation is subsequently supported through Pareto efficiency or sensitivity analysis, a \textbf{two-stage decision-validity failure} arises. In the first stage, scalarisation imposes compensatory trade-offs without establishing that they validly represent actor preferences. In the second stage, Pareto efficiency or sensitivity analysis is used retrospectively to support the resulting optimum, although neither can supply or validate the missing preference foundation. Pareto efficiency does not imply preference validity, and preference validity alone does not establish overall decision validity.

The limitation becomes more fundamental when the fixed directions used in Pareto dominance do not preserve the relevant actor preferences. Conventional Pareto dominance assumes designated directions of improvement, such as ``more is better'' or ``less is better.'' Actors may instead prefer target values, acceptable intervals, saturation levels, or interior optima. Multiple actors may also assign different or opposing preference directions to the same performance. Under these conditions, performance-space dominance need not preserve the ordering induced by aggregated actor preference.

Pareto analysis therefore cannot guarantee that the alternative maximising aggregated actor preference is retained. The resulting Pareto set may exclude the best-fit-for-common-purpose solution altogether. Only when aggregated preference is strictly monotonic in the same directions as the Pareto ordering is a preference-maximising solution guaranteed to be Pareto-efficient.

Pareto analysis may support performance-space ordering and exploration, but ordering non-dominated alternatives is not equivalent to preference-consistent deciding \cite{ozerwolfert2026}. It cannot repair an unvalidated scalarisation, supply missing actor preferences, or establish a decision-valid group optimum. These limitations require a fundamentally different construction of the final decision criterion: system performances must first be transformed into explicitly modelled actor preferences, made commensurate, and only then aggregated.

Section~4 develops this preference-based alternative. The distinction established in this section can therefore be summarised as follows:

\begin{quote}
\textbf{If a scalarisation lacks a preference-valid foundation, Pareto analysis cannot supply or repair it; their combined application therefore cannot guarantee a decision-valid MCDM outcome.}
\end{quote}

\section*{4. From System Performance to Valid Actor Preference}

Sections~2 and~3 showed that neither monetisation nor cost-function scalarisation guarantees decision validity. Monetisation converts heterogeneous performances into quantities of money, while cost-function scalarisation combines transformed performance measures in a weighted numerical criterion. Neither operation by itself establishes valid preference measurement, preference commensurability, actor acceptability, or multi-actor preference aggregation. The constructive alternative is to transform every decision-relevant system performance, including financial cost, into explicitly modelled actor preference before aggregation.

Following Preference Function Modelling (PFM), preference is fundamentally relational and contextual rather than an absolute scalar quantity \cite{barzilai2005,barzilai2010,barzilai2022}. An actor does not possess an absolute quantity such as ``73 units of preference.'' Rather, one performance level or alternative is preferred relative to another. Preference is therefore not intrinsically a scalar quantity, although its relations can be represented numerically on an affine preference scale and used computationally.

\subsection*{4.1 From Performance to Preference, Feasibility, and Acceptability}

For design alternative \(x\), actor \(k\)'s preference concerning performance \(i\) is represented through the performance-to-preference mapping

\begin{equation}
\label{eq:preference-mapping}
p_{k,i}(x)
=
P_{k,i}\bigl(f_i(x)\bigr),
\qquad
p_{k,i}(x)\in[0,100],
\end{equation}

where \(f_i(x)\) is system performance \(i\), \(P_{k,i}\) is the actor-specific performance-to-preference function, and \(p_{k,i}(x)\) is the resulting raw preference representation. Within the specified actor--performance domain and decision context, the relatively least- and most-preferred performance levels are assigned the reference values \(0\) and \(100\), respectively. These endpoints define the chosen origin and unit of an affine preference representation; they do not represent absolute amounts of preference. In particular, \(p_{k,i}(x)=0\) denotes the relatively least-preferred performance level within the specified domain, not zero physical performance or the absence of preference. Likewise, \(p_{k,i}(x)=100\) denotes the relatively most-preferred performance level within that domain and context, not an absolute maximum quantity of preference.

Although the actor--performance-specific preference representations share the numerical range \([0,100]\), this common range does not by itself establish commensurability across actors and performances. Each raw preference representation belongs initially to its own one-dimensional affine preference space and must therefore be transformed into the common Linear Preference Space before valid preference aggregation can take place.

The same performance-to-preference mapping applies to all decision-relevant performances:

\begin{align}
C(x)
&\longrightarrow
P_{k,C}\bigl(C(x)\bigr),
&
T(x)
&\longrightarrow
P_{k,T}\bigl(T(x)\bigr),
\\
E(x)
&\longrightarrow
P_{k,E}\bigl(E(x)\bigr),
&
S(x)
&\longrightarrow
P_{k,S}\bigl(S(x)\bigr).
\end{align}

This mapping synthesises \emph{what the system performs} with \emph{how actors value that performance}. Cost is consequently an object of actor-preference evaluation, not the universal evaluation scale for all other performances.

A decision-valid design--decision model must also distinguish capability, feasibility, desirability, and acceptability. System capability is represented by the performance vector
\begin{equation}
\mathbf{F}(x)
=
\left[
f_1(x),\ldots,f_I(x)
\right],
\end{equation}
while technical, physical, regulatory, and institutional feasibility is imposed through
\begin{equation}
\label{eq:feasibility-condition}
g_f\bigl(\mathbf{F}(x)\bigr)
\leq
0.
\end{equation}

Desirability is represented by the actor-specific preferences \(P_{k,i}(f_i(x))\), whereas acceptability requires
\begin{equation}
\label{eq:acceptability-condition}
P_{k,i}\bigl(f_i(x)\bigr)
\geq
\underline{P}_{k,i}.
\end{equation}

The feasible and acceptable design--decision space is therefore
\begin{equation}
\label{eq:feasible-acceptable-space}
\mathcal{S}_{f,a}
=
\left\{
x\in X:
g_f\bigl(\mathbf{F}(x)\bigr)\leq0
\ \land\
P_{k,i}\bigl(f_i(x)\bigr)\geq\underline{P}_{k,i}
\quad
\forall\,k,i
\right\}.
\end{equation}

Feasibility concerns permissible system behaviour, whereas acceptability concerns what actors are willing to accept. Their separation addresses the compensation problem identified in Section~2: favourable monetary or preference contributions cannot compensate for an impossible or actor-unacceptable outcome. The integration of capability, feasibility, desirability, acceptability, and preference-based solvability is provided by the ODESYS-in-FIVES framework \cite{wolfert2023}.

\subsection*{4.2 Preference Aggregation, Not Performance Scalarisation}

The raw preference representations \(p_{k,i}(x)\) belong to actor--performance-specific one-dimensional affine preference spaces. Although they share the numerical range \([0,100]\), this common range does not make them commensurable across actors and performances. They cannot therefore be added directly merely because they are represented numerically.

Each raw preference representation is transformed into a \(z\)-normalised coordinate in the common Linear Preference Space (LPS):

\begin{equation}
\label{eq:normalised-preference}
z_{k,i}(x)
=
\mathcal{N}_{k,i}
\left(
P_{k,i}\bigl(f_i(x)\bigr)
\right),
\end{equation}

where \(\mathcal{N}_{k,i}\) is the affine normalisation applied to actor \(k\)'s preference concerning performance \(i\). The resulting coordinates \(z_{k,i}(x)\) share a stable zero reference and a commensurate preference unit, thereby establishing the conditions required for preference-consistent aggregation across actors and performances \cite{barzilai2010,wolfert2026unique}.

The \textbf{A-FINE Aggregator}, denoted by \(\mathcal{A}\), combines these commensurate preference coordinates through a linear weighted centroid. The resulting overall aggregated-preference function is

\begin{equation}
\label{eq:aggregated-actor-preference}
Z(x)
=
\mathcal{A}
\left(
\left\{z_{k,i}(x)\right\},
\left\{w'_{k,i}\right\}
\right)
=
\sum_{i=1}^{I}
\sum_{k=1}^{K}
w'_{k,i}\,z_{k,i}(x),
\qquad
w'_{k,i}\geq0,
\qquad
\sum_{i=1}^{I}
\sum_{k=1}^{K}
w'_{k,i}=1.
\end{equation}

The effective contribution weight \(w'_{k,i}=w_k w_{k,i}\) combines the global importance \(w_k\) of actor \(k\) with the local importance \(w_{k,i}\) assigned by that actor to performance \(i\). The weights satisfy \(\sum_{k=1}^{K}w_k=1\) and \(\sum_{i=1}^{I}w_{k,i}=1\) for every actor \(k\), which implies \(\sum_{i=1}^{I}\sum_{k=1}^{K}w'_{k,i}=1\). Thus, \(Z(x)\) is the weighted-centroid representation of aggregated actor preference in the LPS. It is not an absolute quantity of preference; its meaning derives from the underlying preference relations and their common affine representation.

This construction is fundamentally different from the cost-function scalarisation examined in Section~3. Cost-function scalarisation proceeds directly from system performances to a constructed scalar objective, \(\mathbf{F}(x)\longrightarrow J(x)\). Preference aggregation instead follows the sequence \(\mathbf{F}(x)\longrightarrow P_{k,i}\bigl(f_i(x)\bigr)\longrightarrow z_{k,i}(x)\longrightarrow Z(x)\).

The governing distinction is therefore:

\begin{quote}
\textbf{Cost-function scalarisation forms a weighted sum of transformed performance measures; decision-valid MCDM aggregates validly modelled and commensurate actor preferences concerning those performances.}
\end{quote}

he corresponding decision principle is the \textbf{Integrative Maximisation of Aggregated Preferences (IMAP)} \cite{Wolfert2026a}; see also the Data Availability Statement for the corresponding open-source implementation. The best-fit-for-common-purpose alternative is

\begin{equation}
\label{eq:max-aggregated-preference}
x^*
=
\arg\max_{x\in\mathcal{S}_{f,a}}
Z(x).
\end{equation}

IMAP first models actor-specific preferences, represents them in the common LPS, and then applies the A-FINE Aggregator over the feasible and acceptable design--decision space. The distinction from scalarisation therefore does not depend on whether the final computational criterion is scalar-valued, but on what that criterion represents and how it is constructed. Monetisation produces monetary quantities, cost-function scalarisation produces a weighted numerical criterion, and IMAP produces a preference-consistent representation of aggregated actor preference. The alternative \(x^*\) is consequently selected through preference aggregation rather than performance scalarisation.

\section*{5. Conclusion}

This article examined whether monetisation and cost-function scalarisation can guarantee valid multi-criteria decision-making. The central issue is not the use of a numerical criterion or optimisation as such, but whether that criterion validly represents the actor preferences required for deciding.

Section~2 established the first decision-validity failure. Monetisation converts heterogeneous system performances into monetary equivalents, but the resulting quantities remain amounts of money. Financial cost is an object-related system performance, whereas value-for-money is an actor's subjective and contextual evaluation of that cost. Cost-benefit analysis and net present value can provide internally consistent monetary comparisons, while their conversion coefficients encode assumptions about the importance, exchangeability, and compensation of performances. These assumptions are not thereby validated as actor preferences.

The first central conclusion is therefore:

\begin{quote}
\textbf{Monetary commensurability is not preference commensurability: money is not value-for-money.}
\end{quote}

Financial cost is one system performance alongside duration, safety, environmental impact, quality, and reliability. Monetisation elevates money from one performance dimension into the universal representation of all performances, thereby replacing the original multi-performance problem with a single-objective optimisation of an imposed monetary representation. When all terms are actual financial expenditures or revenues, their summation can provide a valid total-cost model. Such a model identifies the least costly feasible system, but it does not, without an explicit preference model, identify the system most preferred by the relevant actors.

ection~3 established the corresponding decision-validity failure of cost-function scalarisation. Replacing heterogeneous performances, including financial cost, with points, normalised scores, losses, penalties, utilities, or fitness measures does not create a preference-neutral model. These transformations determine preference directions, intensities, thresholds, and marginal effects, while their weights determine compensatory trade-offs. A weighted cost function is therefore an implicit preference model, even when it is not monetised. Valid aggregation consequently requires the underlying actor preferences to be explicitly modelled, made commensurate in a common preference domain, and aggregated through a preference-consistent method.

The second central conclusion is therefore:

\begin{quote} \textbf{A weighted scalar cost function embeds subjective preference information but does not constitute valid preference aggregation; decision validity requires preference-consistently aggregated actor preferences.} \end{quote}

Numerical scoring does not by itself establish preference measurement, and normalisation does not establish preference commensurability. Effective trade-offs are jointly determined by score construction and weighting. Different combinations of scales, reference points, normalisations, and weights can consequently produce different optima without any underlying actor preference having changed. A mathematically unique scalar optimum is therefore unique only relative to its numerical construction; it is not necessarily a unique preference optimum or a decision-valid outcome.

Pareto efficiency and sensitivity analysis cannot retrospectively repair an unvalidated scalarisation. Scalarisation first imposes compensatory trade-offs without valid preference measurement; Pareto ordering or parameter variation subsequently evaluates the resulting model without supplying the missing actor preferences. The two-stage failure is therefore cumulative: scalarisation constructs the criterion without a validated preference foundation, and Pareto analysis cannot provide that foundation afterwards.

The resulting decision-validity conclusion is:

\begin{quote}
\textbf{Neither a scalar optimum nor its Pareto efficiency establishes decision validity when the underlying actor preferences have not been validly modelled and aggregated.}
\end{quote}

Section~4 provided the constructive alternative. Decision-valid MCDM must move explicitly from object-related system performances to subject-related actor preferences. Each performance, including financial cost, must first be evaluated through actor-specific performance-to-preference functions. System feasibility must remain distinct from actor acceptability so that favourable performance elsewhere cannot compensate for an infeasible or actor-unacceptable outcome.

The resulting actor--performance-specific preference representations must subsequently be transformed into a common Linear Preference Space with a stable zero reference and a commensurate preference unit. The A-FINE Aggregator combines the resulting \(z\)-normalised preference scores as commensurate coordinates in the Linear Preference Space through a weighted centroid, after which the \textbf{Integrative Maximisation of Aggregated Preferences (IMAP)} optimisation method identifies the feasible and acceptable alternative that maximises aggregate preference. 

\begin{quote}
\textbf{Decision-valid MCDM constitutes preference aggregation, not performance scalarisation.}
\end{quote}

Decision-valid analysis consequently proceeds from system performance to actor preference, then to a common preference domain, valid preference aggregation, and ultimately a best-fit decision. The appropriate final decision criterion is neither the unexamined minimisation of monetised quantities nor the weighted summation of transformed performance scores. It is the maximisation of validly aggregated actor preference over the feasible and acceptable design--decision space. The resulting best-fit-for-common-purpose alternative is determined by what the system can perform, what the actors prefer, what is feasible, and what is mutually acceptable.

\section*{Data Availability}

The contemporary Preferendus software environment comprises Preferendus-IMAP, which uses the A-FINE Aggregator's unique, decision-valid preference-aggregation method to provide preference-consistent design--decision optimisation. The software is publicly available through the following GitHub organisations:

\noindent
\url{https://github.com/Preferendus} and
\url{https://github.com/TUDelft-IMAP-IGS}.


\end{document}